\documentclass[12pt]{amsart}
\usepackage{latexsym}           
\usepackage{amssymb}
\usepackage{epsfig}             
\usepackage{amsfonts}

\newcommand{\IA}{\mathbb{A}}

\newcommand{\IP}{\mathbb{P}}                                     
\newcommand{\IR}{\mathbb{R}}                           
\newcommand{\IC}{\mathbb{C}}
\newcommand{\IZ}{\mathbb{Z}}

\newcommand{\cK}{\mathcal{K}}

\newcommand{\R}{\mathcal{R}}
\newcommand{\A}{\mathcal{A}}
\newcommand{\cB}{\mathcal{B}}
\newcommand{\cC}{\mathcal{C}}
\newcommand{\bcC}{\boldsymbol{\mathcal{C}}}
\newcommand{\cD}{\mathcal{D}}
\newcommand{\cE}{\mathcal{E}}

\newcommand{\cN}{\mathcal{N}}

\newcommand{\G}{\mathcal{G}}

\newcommand{\g}{       \mathfrak{g}     }

\newcommand{\gs}{      \mathfrak{g}^*       }
\newcommand{\gks}{     \mathfrak{g}_k^*   }

\newcommand{\lb}{       \mathfrak{b}            }
\newcommand{\lbk}{      \mathfrak{b}_k          }

\newcommand{\lbks}{     \mathfrak{b}_k^*        }

\newcommand{\lk}{\mathfrak{k}}

\newcommand{\Lia}{{\text{\rm Lie}}}     

\newcommand{\lt}{\mathfrak{t}}

\newcommand{\fus}{\circledast}

\newcommand{\csan}{\overline\Delta_r} 

\newcommand{\bD}{{\bf D}}

\newcommand{\pf}{\begin{bpf}}

\newcommand{\pfms}{\begin{bpfms}}
\newcommand{\epf}{\end{bpf}\hfill$\square$\\}           
\newcommand{\epfms}{\end{bpfms}\hfill$\square$\\}               

\newcommand{\flt}{{\text{flat}}}
\newcommand{\Vect}{{\text{Vect}}}             
\newcommand{\lu}{\mathfrak{u}}
\newcommand{\wt}{\widetilde}

\newcommand{\wh}{\widehat}
\newcommand{\al}{\alpha}
\newcommand{\be}{\beta}
\newcommand{\ga}{\gamma}

\newcommand{\bfe}{{\bf e}}
\newcommand{\bfd}{{\bf d}}

\newcommand{\bs}{{\bf S}}

\newcommand{\ominfo}{{\Omega^1_{C^\infty}}}

\newcommand{\Sect}{\text{\rm Sect}}
\newcommand{\Ssect}{\wh {\text{\rm Sect}}}

\newcommand{\Ad}{\text{\rm Ad}}

\newcommand{\pr}{\text{\rm pr}}         

\newcommand{\res}{{\text{\rm Res}}}
\newcommand{\Res}{{\text{\rm Res}}}

\newcommand{\Hom}{\text{\rm Hom}}

\newcommand{\reg}{{\text{\footnotesize\rm reg}}}

\newcommand {\flb}{\lbrack\!\lbrack}
\newcommand {\frb}{\rbrack\!\rbrack}

\newcommand {\we}{}

\newcommand{\pir}{\pi_{\hbox{\footnotesize\rm res}}}
\newcommand{\pii}{\pi_{\hbox{\footnotesize\rm irr}}}
\newcommand{\spq}{/\!\!/}

\def\mapright#1{\smash{
        \mathop{\longrightarrow}\limits^{#1}}}

\def\mapdown#1{\Big\downarrow
        \rlap{$\vcenter{\hbox{$\scriptstyle#1$}}$}}

\theoremstyle{plain}
\newtheorem {hypo}{\bf\hspace{-\parindent}Hypothesis}
\newtheorem*{plainthm}{Theorem}
\newtheorem {thm}{Theorem}
\newtheorem {prop}[hypo]{Proposition}

\newtheorem {cor}[hypo]{Corollary}
\newtheorem*{plaincor}{Corollary}
\newtheorem {lem}[hypo]{Lemma}

\theoremstyle{definition}
\newtheorem {defn}[hypo]{Definition}
\newtheorem {eg}[hypo]{Example}

\theoremstyle{remark}
\newtheorem {rmk}[hypo]{Remark}

\begin{document}


\title[Quasi-Hamiltonian Geometry of Meromorphic Connections]
{Quasi-Hamiltonian Geometry of \\ Meromorphic Connections}
\author{Philip Boalch}
\address{I.R.M.A., Universit\'e Louis Pasteur et C.N.R.S.\\
7 rue Ren\'e Descartes\\
67084 Strasbourg Cedex\\
France}
\email{boalch@math.u-strasbg.fr}

\begin{abstract}
For each connected complex reductive group G,
we find a family of new examples of complex 
quasi-Hamiltonian G-spaces with G-valued moment maps.
These spaces arise naturally as moduli spaces of 
(suitably framed) meromorphic connections
on principal G-bundles over a disc, and they generalise the 
conjugacy class example
of Alekseev--Malkin--Meinrenken (which appears in the simple 
pole case).
Using the `fusion product' in the theory this gives a
finite dimensional
construction 
of the natural symplectic structures on the spaces
of monodromy/Stokes data  
of meromorphic connections over
arbitrary genus Riemann surfaces, together 
with a new proof of the symplectic nature of isomonodromic deformations
of such connections.
\end{abstract}


\maketitle


\renewcommand{\baselinestretch}{1.0}            
\normalsize

\begin{section}{Introduction}

The quasi-Hamiltonian approach \cite{AMM} to 
constructing symplectic moduli spaces of
flat connections on $G$-bundles over surfaces 
involves ``fusing'' together some basic pieces
and then using a reduction procedure to obtain the symplectic moduli space.
Just two types of such basic quasi-Hamiltonian $G$-spaces are needed
to construct all the moduli spaces considered in \cite{AMM}:
conjugacy classes $\cC\subset G$ and the 
{\em internally fused double} $\bD\cong G\times G$.
Indeed, given a genus $g$ surface $\Sigma$ with $m$ boundary components the
quasi-Hamiltonian reduction 
\begin{equation} \label{eq: first fusion}
(\bD\fus\cdots\fus\bD\fus\cC_1\fus\cdots\fus\cC_m)\spq G
\end{equation}
of the quasi-Hamiltonian fusion of $g$ copies of $\bD$ and $m$ conjugacy
classes $\cC_i$ has a symplectic structure and is isomorphic 
to the moduli space $\Hom_{\bcC}(\pi_1(\Sigma),G)/G$ 
of representations of the fundamental group of $\Sigma$
with holonomy around the $i$th boundary component restricted to lie in $\cC_i$.
Such symplectic moduli spaces have been much studied and in particular 
there are alternative
finite dimensional constructions, cf.  \cite{Kar92,FR,AMR,GHJW}.
A beautiful feature of the approach of \cite{AMM}
is that  the quasi-Hamiltonian moment map condition in the reduction 
\eqref{eq: first fusion} is precisely
the monodromy relation in $\Hom_{\bcC}(\pi_1(\Sigma),G)$.

The aim of this paper is to use the quasi-Hamiltonian approach to give a
finite dimen- sional construction of the natural symplectic structures on more
general moduli spaces where there is currently no other finite dimensional
method. 
This is a continuation of \cite{smid} where the Atiyah--Bott
infinite dimensional approach to moduli spaces of flat connections was
extended to allow singularities in the connections, thereby constructing
symplectic structures on moduli spaces of flat singular connections.
(Such moduli spaces are naturally isomorphic both to spaces of
meromorphic connections with arbitrary order poles over Riemann surfaces, and
to spaces of monodromy/Stokes data, naturally generalising the space of
fundamental group representations above.)

Due to the quasi-Hamiltonian fusion procedure (which, on the level of 
surfaces,
amounts to gluing two surfaces with one boundary component into two of the
holes of a three-holed sphere) we only need to understand moduli
spaces of meromorphic connections on $G$-bundles over a disc having just one
pole. 
This leads to an infinite family of new basic pieces parameterised by the pole
order $k$. These may then be fused together (and with some copies of $\bD$) 
to construct more general moduli spaces.
The main new result of this paper
is as follows and will be proved algebraically without referring to 
meromorphic connections or the Stokes phenomenon.
Let $k\ge 2$ be an integer and suppose $U_\pm$ are the full unipotent
subgroups of a pair of Borels $B_\pm\subset G$ intersecting in a maximal torus 
$T=B_+\cap B_-$ of $G$. 

\begin{plainthm}
The manifold $\wt\cC = G\times (U_+\times U_-)^{k-1}\times\lt$
is a complex quasi-Hamiltonian $G\times T$-space with action
$$(g,t)\cdot(C,S_1,\ldots,S_{2k-2},\Lambda) = 
(tCg^{-1},tS_1t^{-1},\ldots,tS_{2k-2}t^{-1},\Lambda)\in\wt\cC$$
(where $S_{\text{odd/even}}\in U_{+/-}, (g,t)\in G\times T$),
moment map $(\mu,e^{-2\pi i\Lambda}):\wt\cC\to G\times T$ where
$$\mu : \wt \cC\to G; \quad 
(C,\bs,\Lambda)\mapsto C^{-1}S_{2k-2}\cdots S_2S_1e^{2\pi i\Lambda}C=D^{-1}E,
$$ 
and two-form
\begin{equation} \label{eq: first omega exprn}
\boxed{\ \ \omega = \frac{1}{2}(\overline\cD,\overline\cE) + 
\frac{1}{2}\sum_{j=1}^{k-1} (\cD_j,\cD_{j-1}) - (\cE_j,\cE_{j-1}) \ \ }
\end{equation}
\noindent
where $\overline\cD=D^*\overline\theta, 
\overline\cE=E^*\overline\theta,
\cD_j=D_j^*\theta,\cE_j=E_j^*\theta\in\Omega^1(\wt \cC,\g)$
for maps $D_j,E_j:\wt \cC\to G$ defined by
$D_j=\epsilon^{-j}S_{2k-1-j}^{-1}\cdots S^{-1}_{2k-3}S^{-1}_{2k-2}C,
E_j=\epsilon^{j+2-2k}S_j\cdots S_2S_1\epsilon^{2k-2}C,$
$D:=D_{k-1},E:=E_{k-1}, E_0=D_0:=C$ where 
$\epsilon:= e^{\frac{\pi i \Lambda}{k-1}}$.
Moreover for each choice of $\Lambda\in\lt$ the reduction
$$\cC :=  (\wt\cC\vert_\Lambda)/ T \cong (G\times(U_+\times U_-)^{k-1})/T$$
is a  complex quasi-Hamiltonian $G$-space.
\end{plainthm}

(In the body of the paper a more symmetrical, $\epsilon$-free, notation will
be used.)
Thus for each pole order there are new quasi-Hamiltonian spaces 
$\cC$ and $\wt \cC$, 
the first of which arises simply from the second upon reducing
by a torus, and depends on a choice of $\Lambda$. 


For example, in the order two pole case $k=2$, if we define
$b_- = e^{-\pi i \Lambda} S_2^{-1},$ and 
$b_+ = e^{-\pi i \Lambda} S_1e^{2\pi i\Lambda}
$ then 
$$\wt\cC \cong G\times G^*,\qquad\mu=C^{-1}b_-^{-1}b_+C,\qquad
\omega = \frac{1}{2}(\overline\cD,\overline\cE) + 
\frac{1}{2}(\cD,\gamma) - \frac{1}{2}(\cE,\gamma)$$
where $G^*$ is the simply connected Poisson Lie group dual to $G$
and 
$D=b_-C, E=b_+C, \gamma=C^*\theta$.
In general 
the quotient $\wt\cC/G$ has an induced Poisson structure \cite{AKM} and for
$k=2$ this coincides with 
standard Poisson structure on $G^*$.
Also we will see that for $k=2$ the additive analogue $\wt O$ of $\wt \cC$ is
the cotangent bundle $T^*G$.

To understand the geometrical origins of these spaces we remark that the
Stokes multipliers $\{S_i\}$ arise from the irregular Riemann-Hilbert
correspondence; in Theorem \ref{thm: what c is} below 
it is explained how $\wt\cC$ is isomorphic to a
moduli space of (framed) meromorphic connections with fixed irregular
type on $G$-bundles over a disc.

Given a genus $g$ compact Riemann surface $\Sigma$ with a divisor
$D=\sum_{i=1}^mk_i(a_i)$ having each $k_i\ge 1$ the above theorem enables
one to construct
symplectic moduli spaces of monodromy data for meromorphic connections
on $\Sigma$ of the form
\begin{equation} \label{eq: second fusion}
(\bD\fus\cdots\fus\bD\fus\wt\cC_1\fus\cdots\fus\wt\cC_m)\spq G
\end{equation}
with $g$ factors of $\bD$,
as constructed in \cite{smid} from an infinite dimensional viewpoint.
Summing the quasi-Hamiltonian two-forms on each factor in 
\eqref{eq: second fusion} together with the fusion terms gives an explicit
expression for the symplectic structure on the manifold 
\eqref{eq: second fusion}.
Such an expression has also been obtained directly in the recent preprint
\cite{Krich-imds}, however the approach here gives an algebraic
proof that it is indeed a symplectic structure.

In section \ref{sn: additive} we will recall the additive analogues $O,\wt O$
of the spaces $\cC,\wt\cC$ for each $k$. These are symplectic manifolds of 
matching dimensions ($\dim\cC=\dim O, \dim\wt O= \dim\wt \cC$).
Indeed $O$ is just a generic coadjoint orbit of the group 
$G_k:=G(\IC[z]/z^k)$ of $(k-1)$-jets of bundle automorphisms
and so this nicely extends the idea that the conjugacy classes are the
multiplicative analogue of coadjoint orbits of $G$.
The {\em extended orbits} $\wt O$ are larger by $2\dim T$ and give rise to the
orbits $O$ upon performing  a symplectic quotient by $T$.

In the genus zero case the spaces $O,\wt O$ enable one to construct global
symplectic moduli spaces of meromorphic connections on trivial $G$-bundles
as symplectic quotients 
of the form $(\wt O_1\times\cdots\times \wt O_m)\spq G$, and in fact such
moduli spaces fill out a dense subset of a component of the full moduli space.
The main result of \cite{smid} then leads immediately to:

\begin{plaincor}
The (global) irregular Riemann-Hilbert map 
\begin{equation} \label{eq: irhc}
\nu : (\wt O_1\times\cdots\times \wt O_m)\spq G\hookrightarrow
(\wt\cC_1\fus\cdots\fus\wt\cC_m)\spq G
\end{equation}
associating
monodromy/Stokes data to a meromorphic connection on a trivial $G$-bundle
over $\IP^1$ is a symplectic map (provided the symplectic structure on the
right-hand side is divided by $2\pi i$).
Moreover both sides are naturally Hamiltonian $T^{\times m}$-spaces and $\nu$
intertwines these actions and their moment maps.
\end{plaincor}

This map $\nu$ depends heavily on the chosen pole positions $a_i$
(and on a choice of `irregular type' at each pole). 
However from the formula \eqref{eq: first omega exprn} the symplectic
structure on the right-hand side of \eqref{eq: irhc} is manifestly independent
of these choices.
This shows that the isomonodromy connection is a symplectic connection, as was
shown in \cite{smid} from a de\! Rham point of view.

For example the case with two poles of order two on
$\IP^1$ looks as follows. 
Since $\wt O\cong T^*G$ for $k=2$, the left-hand side of 
\eqref{eq: irhc} is also isomorphic to the cotangent bundle $T^*G$.
The right-hand side of \eqref{eq: irhc} turns out 
(Proposition \ref{prop: sym dgpoid}) to be isomorphic as a symplectic
manifold to the symplectic double groupoid $\Gamma$ of $G$ and $G^*$, described
in \cite{LuW-sdgpoid}. 
Thus in this case $\nu$ is a (transcendental) embedding
$T^*G\hookrightarrow
\Gamma$ between these (essentially algebraic) symplectic manifolds for each
choice of pole positions and irregular types.



\renewcommand{\baselinestretch}{1}		
\small

\ 

{\em Acknowledgements.}
The author would like 
to thank Anton Alekseev and Ping Xu for the opportunity to talk about the
$k=2$ version of these results at the conference on Poisson Geometry,
E.S.I. Vienna, June 2001.
The author is supported by the
European Differential Geometry Research Training Network (EDGE)
HPRN-CT-2000-00101.

\renewcommand{\baselinestretch}{1}		
\normalsize

\end{section}

\begin{section}*{Notation/Conventions}
Throughout this paper 
$G$ is a connected complex reductive group with maximal torus $T$ and
corresponding Lie algebras $\lt\subset\g$.
$B_\pm$ denote a pair of opposite Borels subgroups 
with $B_+\cap B_- = T$ 
and $U_\pm\subset B_\pm$ denote their full unipotent subgroups,
with corresponding Lie algebras $\lu_\pm\subset\lb_\pm$.

Let $(,):\g\otimes\g\to \IC$ be a symmetric nondegenerate invariant
bilinear form. 
(Note that invariance implies $(,)$ restricts to zero on
$\lu_\pm\otimes\lu_\pm$ and to a pairing on each of
$\lu_+\otimes\lu_-, \lu_-\otimes\lu_+, \lt\otimes\lt$.)

$\theta,\overline\theta\in\Omega^1(G,\g)$ denote the tautological 
left and right invariant $\g$-valued holomorphic one-forms on $G$
respectively (so in any representation 
$\theta=g^{-1}dg, \overline\theta=(dg)g^{-1}$).

Generally if $\A,\cB,\cC\in\Omega^1(M,\g)$ are $\g$-valued
holomorphic one-forms on a complex manifold $M$ then
$(\A,\cB)\in\Omega^2(M)$ and 
$[\A,\cB]\in\Omega^2(M,\g)$ are defined
by wedging the form parts and pairing/bracketing the Lie algebra parts
(so e.g. $(A\al,B\be)= (A,B)\al\wedge\be$ for
$A,B\in \g, \al,\be\in\Omega^1(M)$).

Define $\A\cB := \frac{1}{2}[\A,\cB]\in \Omega^2(M,\g)$ (which works out
correctly in
any representation of $G$  using matrix multiplication).
Then one has 
$d\theta=-\theta^2, d\overline\theta = \overline\theta^2$.

Define $(\A\cB\cC) = (\A, \cB\cC)\in\Omega^3(M)$ (which is
totally symmetric in $\A, \cB, \cC$). The canonical 
bi-invariant three-form on $G$ is then
$\eta:= \frac{1}{6}(\theta^3)$.

The adjoint action of $G$ on $\g$ will be 
denoted $gXg^{-1}:=\Ad_gX$ for any
$X\in\g, g\in G$.

If $G$ acts on $M$, the fundamental vector field of $X\in\g$ is minus
the tangent to the flow 
$(v_X)_m = -\frac{d}{dt} (e^{Xt}\cdot m)\bigl\vert_{t=0}$, so that the map
$\g\to\Vect_M; X\to v_X$ is a Lie algebra homomorphism. 
(This sign convention differs from \cite{AMM} (and agrees with \cite{AKM}); 
this leads to sign changes in the quasi-Hamiltonian axioms and
the fusion and equivalence theorems.)

\end{section}

\begin{section}{Quasi-Hamiltonian G-spaces}

\begin{defn}[cf. \cite{AMM,AKM}]
A complex manifold $M$ is a 
{\em complex quasi-Hamiltonian $G$-space}
if there is an action of $G$ on $M$, 
a $G$-equivariant map $\mu:M\to G$ (where $G$ acts on itself by
conjugation) and a $G$-invariant holomorphic two-form
$\omega\in \Omega^2(M)$ such that

\noindent(QH1). 
The exterior derivative of $\omega$ is the pullback of the canonical
three-form on $G$:
$$d\omega = \mu^*(\eta).$$

\noindent(QH2).
For all $X\in \g$
$$\omega(v_X,\cdot\,) = \frac{1}{2}\mu^*(\theta+\overline\theta, X)
\in \Omega^1(M).$$

\noindent(QH3).
At each point $m\in M$ the kernel of $\omega$ is
$$\ker\omega_m = \{(v_X)_m\ \bigl\vert\ 
X\in\g \text{ satisfies } gXg^{-1} = -X \text{ where } 
g:=\mu(m)\in G\}.$$
\end{defn}

These axioms are perhaps best motivated in terms of Hamiltonian loop group
manifolds \cite{AMM, MW-Cob}, as we will sketch in Section \ref{sn: derivn}.

A simple but important observation is that if $G$ is abelian 
(and in particular if $G=\{1\}$ is trivial) then these axioms imply 
 that the two-form $\omega$ is a symplectic form.


\begin{eg}[Conjugacy classes \cite{AMM}]
Let $\cC\subset G$ be a conjugacy class, with the conjugation action
of $G$ and moment map $\mu$ given by the inclusion map.
Then $\cC$ is a quasi-Hamiltonian $G$-space with two-form
$\omega$ determined by 
$$\omega_g(v_X,v_Y) = 
\frac{1}{2}\bigl((X,gYg^{-1}) - (Y,gXg^{-1})\bigr) $$
for any $X,Y\in \g, g\in \cC$.
For later use we note that if $g\in \cC$ is fixed and 
we define the surjective map 
$\pi:G \to \cC: C\mapsto C^{-1}gC$  
then one may calculate 
\begin{equation} \label{eq k=1 pullback}
\pi^*(\omega) = \frac{1}{2}(\overline\theta, g\overline\theta
g^{-1})\in\Omega^2(G).
\end{equation}
\end{eg}

\begin{eg}[Internally Fused Double \cite{AMM}]
The space $\bD= G\times G$ is a quasi-Hamiltonian $G$-space 
with $G$ acting by diagonal conjugation 
($g(a,b)=(gag^{-1}, gbg^{-1})$), moment map given by the group commutator
$$\mu(a,b)= aba^{-1}b^{-1}$$
and two-form
$$\omega_\bD=
-\frac{1}{2}(a^*\theta,b^*\overline\theta)
-\frac{1}{2}(b^*\theta,a^*\overline\theta)
-\frac{1}{2}((ab)^*\theta,(a^{-1}b^{-1})^*\overline\theta).$$
\end{eg}

Now let us recall the quasi-Hamiltonian reduction theorem:

\begin{thm}[\cite{AMM}]
Let $M$ be a quasi-Hamiltonian $G\times H$-space with moment map
$(\mu,\mu_H):M\to G\times H$ and suppose
that the quotient by $G$ of the inverse image $\mu^{-1}(1)$  
of the identity under the first moment map is a manifold.
Then the restriction of the two-form $\omega$ to $\mu^{-1}(1)$ 
descends to the {\em reduced space}
\begin{equation*}
M\spq G := \mu^{-1}(1)/G 
\end{equation*}
and makes it into a quasi-Hamiltonian $H$-space. 
In particular, if $H$ is abelian (or in particular trivial) 
then $M\spq G$ is a symplectic manifold. 
\end{thm}


The fusion product, which puts a ring structure on the  
category quasi-Hamiltonian $G$-spaces, is defined as follows. 
(Also reduction at different values of the moment map may be  facilitated
by first fusing with a conjugacy class, analogously to the Hamiltonian case.)

\begin{thm}[\cite{AMM}]\label{thm: fusion}
Let $M$ be a quasi-Hamiltonian $G\times G\times H$-space, 
with moment map $\mu=(\mu_1,\mu_2,\mu_3)$. 
Let $G\times H$ act by 
the diagonal embedding $(g,h)\to (g,g,h)$. 
Then $M$ with two-form 
\begin{equation} \label{eqn: fusion 2form}
\wt{\omega}= \omega - \frac{1}{2}(\mu_1^* \theta, \mu_2^* \overline \theta)
\end{equation}
and moment map
$$\wt{\mu} = (\mu_1\cdot \mu_2,\mu_3):M\to G\times H$$
is a quasi-Hamiltonian $G\times H$-space. 

\end{thm}

We will refer to the extra term subtracted off in \eqref{eqn: fusion 2form} as
the ``fusion term''. 
If $M_i$ is a quasi-Hamiltonian $G\times H_i$ space for $i=1,2$ their 
fusion product $$M_1\fus M_2$$ 
is defined to be the
quasi-Hamiltonian $G\times H_1\times H_2$-space 
obtained from the
quasi-Hamiltonian $G\times G \times H_1 \times H_2$-space 
$M_1\times M_2$ by
fusing the two factors of $G$.


\end{section}

\begin{section}{New Examples}

In this section we will describe the family of quasi-Hamiltonian spaces
$\cC,\wt \cC$
and prove directly that they are such.
The motivation for, and geometrical origins of, 
these spaces will only become clear
in Section \ref{sn: derivn} however, where their infinite dimensional
counterparts are described.

Our main objects of study are the family of complex manifolds
\begin{equation}
\wt \cC := \{
(C,\bfd,\bfe,\Lambda)\in 
G\times(B_-\times B_+)^{k-1}\times\lt \ \bigl\vert\ 
\delta(d_j)^{-1} = e^{\frac{\pi i\Lambda}{k-1}} 
= \delta(e_j) \ \text{for all $j$}\},
\end{equation}
parameterised by an integer $k\ge 2$,
where $\bfd=(d_1,\ldots,d_{k-1}),\bfe=(e_1,\ldots,e_{k-1})$ with
$d_{\text{even}} , e_{\text{odd}}\in B_+$ and 
$d_{\text{odd}},e_{\text{even}}\in B_-$
and where $\delta:B_\pm\to T$ is the homomorphism with kernel $U_\pm$.
This space $\wt\cC$ 
is isomorphic to $G\times(U_+\times U_-)^{k-1}\times\lt$ but 
it will be more convenient to use the above definition.
For the record, in terms of the Stokes multipliers:
$
d_j = \epsilon^{-j}S^{-1}_{2k-1-j}\epsilon^{j-1},
e_j=\epsilon^{j+2-2k}S_j\epsilon^{2k-1-j}$ 
where
$\epsilon:= e^{\frac{\pi i \Lambda}{k-1}}.$

In this description the action of $G$ on $\wt \cC$ given by 
$$g\cdot(C,\bfd,\bfe,\Lambda) = (Cg^{-1},\bfd,\bfe,\Lambda)$$
for $g\in G$, and the action of $T$ is given by
$$t\cdot(C,\bfd,\bfe,\Lambda) = 
(t C,td_1t^{-1},\ldots,td_{k-1}t^{-1},te_1t^{-1},\ldots,te_{k-1}t^{-1}
,\Lambda)$$
for $t\in T$.
Independently these actions are both free, although the combined $G\times T$
action is not.
The maps $D_i,E_i,\mu : \wt \cC\to G$ are now defined as
$$D_i( C,\bfd,\bfe,\Lambda) = d_i\cdots d_1C\qquad (i=0,1,\ldots,k-1)$$
$$E_i( C,\bfd,\bfe,\Lambda) = e_i\cdots e_1C\qquad (i=0,1,\ldots,k-1)$$
$$\mu( C,\bfd,\bfe,\Lambda)= 
C^{-1}d_1^{-1}\cdots d_{k-1}^{-1}e_{k-1}\cdots e_1C.$$
To lighten the notation we will write $D=D_{k-1}, E= E_{k-1}$,
so in particular $\mu=D^{-1}E$.
The main result of this section is:

\begin{thm} \label{mainthm}
The manifold $\wt \cC$ is a quasi-Hamiltonian $G\times T$-space with 
the above action, moment map 
$(\mu,e^{-2\pi i\Lambda}):\wt \cC\to G\times T$ and
two-form: 
\begin{equation} \label{eq: 2form}
\omega = \frac{1}{2}(\overline\cD,\overline\cE) + 
\frac{1}{2}\sum_{i=1}^{k-1} (\cD_i,\cD_{i-1}) - (\cE_i,\cE_{i-1}),
\end{equation}
where
$\overline\cD=D^*(\overline\theta), 
\overline\cE=E^*(\overline\theta),
\cD_i=D_i^*(\theta),\cE_i=E_i^*(\theta)\in\Omega^1(\wt \cC,\g)$.
\end{thm}

In particular, since $T$ is abelian, this implies $\wt \cC$ is a
quasi-Hamiltonian $G$-space with moment map $\mu$ and the same 
two-form.

\begin{rmk}
Observe that $\omega$ is invariant under translations of $\Lambda$
by the lattice 
$L:=\ker(\exp(2\pi i\,\cdot\,):\lt\to T)$.
The quotient $\wt \cC/L\cong G\times(U_+\times U_-)^{k-1}\times T$ is
then an {\em algebraic} quasi-Hamiltonian $G\times T$-space.
Indeed all the formulae above make sense directly on the subvariety of 
$G\times(B_-\times B_+)^{k-1}$ cut out by the
equations $\delta(d_i)\cdot\delta(e_j)=1$, and this subvariety is a
finite covering of $\wt \cC/L$ (corresponding to replacing 
$\exp(\frac{\pi i\Lambda}{k-1})$ by $\exp(2\pi i\Lambda)$).
However we prefer to keep in the choice of $\Lambda$ in order to obtain 
(genuine) Hamiltonian $T$-spaces upon reducing by $G$.
\end{rmk}

Also we have:
\begin{cor}
Suppose a value of $\Lambda$ is fixed. 
Then the reduction 
$$\cC := (\wt\cC\vert_\Lambda)/ T$$
is a complex (algebraic) quasi-Hamiltonian $G$-space. 
\end{cor}
\pf
$\cC$ may also be described as
the quasi-Hamiltonian reduction $(\wt\cC/L)\spq T$ of $\wt\cC/L$
at the value
$\exp(-2\pi i\Lambda)$ of the moment map for the $T$ action.
\epf


Before proving Theorem \ref{mainthm} let us describe 
the order two pole case $k=2$ and the specialisation  to the simple pole
case $k=1$.
If $k=2$ and we define
$b_- = e^{-\pi i \Lambda} S_2^{-1},$ and 
$b_+ = e^{-\pi i \Lambda} S_1e^{2\pi i\Lambda}
$ (so that $b_-^{-1}b_+ = S_2S_1e^{2\pi i\Lambda}$)
then 
\begin{equation} \label{eq: k=2 formulae}
\wt\cC \cong G\times G^*,\qquad\mu=C^{-1}b_-^{-1}b_+C,\qquad
\omega = \frac{1}{2}(\overline\cD,\overline\cE) + 
\frac{1}{2}(\cD,\gamma) - \frac{1}{2}(\cE,\gamma)
\end{equation}
where $D=b_-C, E=b_+C, \gamma=C^*\theta$ and
$$G^* := \{ (b_-,b_+,\Lambda)\in B_-\times B_+\times\lt \ \bigl\vert \  
\delta(b_-)\delta(b_+)=1, \delta(b_+)=\exp(\pi i \Lambda) \}$$
is the Poisson Lie group dual to $G$ (cf. e.g. \cite{smpgafm},
\cite{bafi} Appendix B).

Considering two poles of order two on $\IP^1$ leads to the following
statement, which gives a relationship 
between symplectic double groupoids and meromorphic connections.
\begin{prop} \label{prop: sym dgpoid}
Let $\wt \cC_1$ and $\wt \cC_2$ be two copies of $\wt\cC$ with $k=2$.
Then the quasi-Hamiltonian reduction of the fusion of $\wt\cC_1$ and $\wt\cC_2$
is isomorphic as a symplectic manifold to the symplectic double
groupoid $\Gamma$ of $G$ and $G^*$ appearing in \cite{LuW-sdgpoid}:
\begin{equation} \label{eq: spd}
(\wt\cC_1\fus\wt\cC_2)\spq G\cong \Gamma.
\end{equation}
\end{prop}
\pf
First recall $\wt \cC_i\cong G\times G^*$ as manifolds.
We will assume that the Borels chosen at the first pole 
are opposite to those chosen at the second
(which we may since isomonodromy will give symplectic isomorphisms
with the spaces arising from any other choice of Borels intersecting in $T$).
Thus 
$\wt\cC_1 = \{(C_1,b_-,b_+, \Lambda_1)\vert 
\delta(b_\pm)=e^{\pm\pi i \Lambda_1}\}$
and $\wt\cC_2 = \{(C_2,c_+,c_-, \Lambda_2)\vert 
\delta(c_\pm)=e^{\mp\pi i \Lambda_2}\}$
with $b_\pm,c_\pm\in B_\pm$.
The moment map
on $\wt\cC_1\fus\wt\cC_2$ is
$\mu= C_1^{-1}b_-^{-1}b_+C_1C_2^{-1}c^{-1}_+c_-C_2$.
Writing $h:=C_2C^{-1}_1$ the condition $\mu=1$ becomes
$hb_-^{-1}b_+h^{-1}c_+^{-1}c_-= 1,$ and if
we define $g:=c_-hb_-^{-1}$ then this condition is clearly 
equivalent to $c_+h=gb_+$.
Thus (omitting the $\Lambda$ terms to simplify notation)
we have defined a surjective map 
$$\mu^{-1}(1) \to \Gamma:= 
\{(g,b_-,b_+,h,c_+,c_-)\vert c_\pm h=g b_\pm\}\subset
(G\times G^*)^2$$
whose fibres are precisely the $G$ orbits.
This is the definition of the manifold $\Gamma$ given in
\cite{LuW-sdgpoid}.
The symplectic structures may be shown to agree as follows.

The map $\Gamma\to G\times G; (g,b_-,b_+,h,c_+,c_-)\mapsto 
(gb_-,gb_+)$  expresses $\Gamma$ as the covering of a dense subset
of $G\times G$. This subset is the big symplectic leaf of
the Heisenberg double Poisson structure on $G\times G$ 
and the symplectic structure on $\Gamma$ is defined to be the pullback
of the symplectic structure on this leaf.
An explicit formula for this pullback (i.e. for the symplectic
structure on $\Gamma$) is given in Theorem 3 of 
\cite{AlekMalk94}.
On the other hand we have an explicit formula for the symplectic structure on 
$(\wt\cC_1\fus\wt\cC_2)\spq G$ (involving seven terms, the fusion term
plus three terms \eqref{eq: k=2 formulae} for each factor). 
A straightforward calculation shows these explicit formulae on each
side agree.
\epf

In the simple pole case $k=1$ we {\em define} $\wt \cC = G\times\lt_1$ where 
$\lt_1\subset\lt$ is the complement of the affine root hyperplanes:
$\lt_1:=
\{\Lambda\in\lt\ \bigl\vert \ \al(\Lambda)\notin\IZ \text{ for all
roots } \al \}$.
The correct specialisation of \eqref{eq: first omega exprn} to this case is 
$$\omega = \frac{1}{2}(\overline\cD,\overline\cE) + 
\frac{1}{2}(\cD,\gamma) - \frac{1}{2}(\cE,\gamma)
=2\pi i(\overline\ga,d\Lambda)+ 
\frac{1}{2}(\overline\ga,e^{2\pi i\Lambda}\overline\ga e^{-2\pi i\Lambda})$$
where $D=e^{-\pi i\Lambda}C, E=e^{\pi i\Lambda}C, 
\overline\ga = C^*\overline\theta$.
This is the restriction of the two-form  \eqref{eq: first omega exprn}
to the submanifold  with $d_i,e_i\in T, \Lambda\in\lt_1$ (for any $k$) and 
makes $\wt \cC$ into a quasi-Hamiltonian $G\times T$-space
with moment map $(D^{-1}E,e^{-2\pi i\Lambda})$.
(It also arises as a cross-section of the double in \cite{AMM}.)
Given a fixed choice of $\Lambda$, 
the restriction of $\omega$  to $G\times\{\Lambda\}$ clearly agrees with
\eqref{eq k=1 pullback} and so we deduce
the reduction $\cC :=  (\wt\cC\vert_\Lambda)/ T \cong G/T$
is 
isomorphic as a quasi-Hamiltonian $G$-space to the conjugacy class
through $e^{2\pi i\Lambda}$.

\pfms {\em(of Theorem \ref{mainthm}).}\ \ 
To establish (QH1), that $\mu^*(\theta^3) = 6d\omega$, we observe that the
expression \eqref{eq: 2form} defines a two-form 
on $G\times (B_+\times B_-)^{k-1}$, and working 
algebraically we will view (QH1) as a statement about the
differential algebra generated by the symbols $C,d_i,e_j$, 
using only the restriction that $(d_i^*\theta^3)=(e_j^*\theta^3)=0$
(which follow from the fact that $d_i,e_j$ live in Borel subgroups).
By restriction the result for 
$\wt \cC$ will then follow.
This viewpoint enables us to use induction on $k$.
From the definition one finds 
$\mu^*(\theta^3) = ((\overline \cE - \overline \cD)^3)$ which expands to give
$$ \mu^*(\theta^3)= 
3 (\overline \cD\overline\cD\overline\cE)-
3 (\overline \cD\overline\cE\overline\cE)
+ (\overline\cE^3) - (\overline\cD^3).
$$
On the other hand
$$
2d\omega = 
(\overline \cD\overline\cD\overline\cE)-
(\overline \cD\overline\cE\overline\cE)+ F_{k-1}$$
where
$$
F_{k-1} := \sum_{i=1}^{k-1} 
(\cD_i\cD_{i-1}\cD_{i-1})-
(\cD_i\cD_{i  }\cD_{i-1})-
(\cE_i\cE_{i-1}\cE_{i-1})+
(\cE_i\cE_{i  }\cE_{i-1}),
$$
so that what we must prove is 
$(\overline\cE^3) - (\overline\cD^3) = 3 F_{k-1}$
or equivalently (assuming
$(\overline\cE_{k-2}^3) - (\overline\cD_{k-2}^3) = 3 F_{k-2}$
inductively) that
$(\overline\cE^3) - (\overline\cD^3)$ equals 
\begin{equation} \label{eqn: 1st 3}
(\overline\cE_{k-2}^3) - (\overline\cD_{k-2}^3)
+3\bigl(
(\cD\cD_{k-2}\cD_{k-2})-
(\cD\cD\cD_{k-2})-
(\cE\cE_{k-2}\cE_{k-2})+
(\cE\cE\cE_{k-2})
\bigr).
\end{equation}
To establish this, write $E=b_+E_{k-2}, D=b_- D_{k-2}$ where
$b_+:=e_{k-1}, b_-:= d_{k-1}$. 
(Note we do not necessarily have $b_\pm\in B_\pm$, only that they are in
opposite Borels.)
Thus 
\begin{equation}\label{cEcD}
\cE = E_{k-2}^{-1}\theta_+E_{k-2} + \cE_{k-2},\qquad
\cD = D_{k-2}^{-1}\theta_-D_{k-2} + \cD_{k-2},
\end{equation}
$$\overline\cE = \overline \theta_+  +b_+\overline\cE_{k-2}b_+^{-1},\qquad
\overline\cD = \overline \theta_-  +b_-\overline\cD_{k-2}b_-^{-1},$$
where 
$\theta_\pm = b_\pm^*(\theta), \overline\theta_\pm = b_\pm^*(\overline\theta)$
and so
$$(\overline\cE^3) - (\overline\cD^3) =
 ((\theta_+ + \overline\cE_{k-2})^3) - 
 ((\theta_- + \overline\cD_{k-2})^3) = $$
\begin{equation} \label{eqn: 2nd 3}
(\overline\cE_{k-2}^3) - (\overline\cD_{k-2}^3)
+3\bigl(
(\theta_+\theta_+\overline\cE_{k-2})+
(\theta_+\overline\cE_{k-2}\overline\cE_{k-2})-
(\theta_-\theta_-\overline\cD_{k-2})-
(\theta_-\overline\cD_{k-2}\overline\cD_{k-2})
\bigr)
\end{equation}
using the fact that $(\theta_\pm^3)=0$.
Thus we must show that the coefficients of $3$ in 
\eqref{eqn: 1st 3} and \eqref{eqn: 2nd 3} are the same;
this however follows easily by substituting the expressions 
\eqref{cEcD} for 
$\cE, \cD$ into $\eqref{eqn: 1st 3}$ and expanding.
Finally the $k=2$ case may be proved directly, justifying the
induction; namely we must show 
$(\overline\cE^3) - (\overline\cD^3) = 3 F_{1}$
and this comes about simply by expanding both sides in terms of $b_\pm$ and
$C$. (The $k=1$ case is similar.) 


Next we will check (QH2) for the $G$-action.
Choose $X\in\g$ and an arbitrary holomorphic vector field $Y$ on
$\wt \cC$. We will denote derivatives along $v_X$ by primes and along
$Y$ by dots, so e.g. 
$\dot\cD_i = \langle Y, \cD_i \rangle \in \g$ and
$\cE'_j = \langle v_X, \cE_j\rangle \in \g$
(and in any representation of $G$ we have
$\dot\cD_i = D_i^{-1}\dot D_i$ etc).
By definition of the action 
$\cD'_i=\cE'_i=X$ for all $i$, and 
$\overline\cD'=DXD^{-1}, \overline\cE'=EXE^{-1}$.
Thus 
$$2\omega(v_X,Y) = 
(DXD^{-1},\dot{\overline\cE}) - (EXE^{-1},\dot{\overline\cD})
+\sum_{i=1}^{k-1}\Bigl(
X,\dot\cD_{i-1}-\dot\cD_i-\dot\cE_{i-1}+\dot\cE_i
\Bigr)$$
which simplifies to 
$(X,D^{-1}\dot{\overline\cE}D-E^{-1}\dot{\overline\cD}E  
-\dot\cD + \dot\cE)$.
On the other hand, since $\mu=D^{-1}E$:
$$\langle(\mu^*\theta+\mu^*\overline\theta,X),Y\rangle=
(\mu^{-1}\dot\mu+\dot\mu\mu^{-1},X) =
(\dot\cE-\dot\cD+D^{-1}\dot{\overline\cE} D-E^{-1}\dot{\overline\cD} E,X)$$
so we have established (QH2) for the $G$-action.

For the $T$-action, if $X\in\lt$  then the derivatives along the 
corresponding fundamental
vector field $v_X$ (for the $T$ action) are:
$\dot{\overline\cD}_i=\dot{\overline\cE}_i = -X ,
\dot\cD_i=-D_i^{-1}XD_i, \dot\cE_i=-E_i^{-1}XE_i.$
Thus for any vector field $Y$ on $\wt \cC$
$$2\omega(v_X,Y)=
\bigl(X,-\overline\cE'+ \overline\cD'
+\sum_{i=1}^{k-1}
-D_i\cD_{i-1}'D_i^{-1}+E_i\cE'_{i-1}E_i^{-1}
+D_{i-1}\cD_{i}'D_{i-1}^{-1}-E_{i-1}\cE'_{i}E_{i-1}^{-1}\bigr)$$
where the primes denote the derivatives along $Y$.
Now $D_i=d_iD_{i-1}$ so that 
$D_i\cD_{i-1}'D_i^{-1}=\overline\cD'_i-\overline\delta'_i$ and
$D_{i-1}\cD_{i}'D_{i-1}^{-1}=\overline\cD'_{i-1}+\delta'_i$ 
(and similarly for the $\cE_i$'s), where 
$\delta_i:=d_i^*\theta$ etc.
Substituting thus shows
$$2\omega(v_X,Y)= \bigl(X,
\sum_{i=1}^{k-1}
\delta'_i+\overline\delta'_i-\varepsilon'_i-\overline\varepsilon'_i\bigr).$$
Since $X\in\lt$ we may take the $\lt$ component of the right-hand side 
yielding 
$$\omega(v_X,Y)= -(2\pi i)(X,\Lambda')=-(2\pi
i)\langle(d\Lambda,X),Y\rangle$$
which is what appears on the right-hand side of (QH2) if the 
 moment map is $e^{-2\pi i\Lambda}$. 
 
The proof of the minimal degeneracy condition (QH3) is rather complicated 
so has been put in the appendix.
\epfms
\end{section}

\begin{section}{Derivation} \label{sn: derivn}

In this section we will explain how the quasi-Hamiltonian spaces
$\cC,\wt\cC$ were found.
In brief the extension of the Atiyah--Bott symplectic structure to the
meromorphic case in \cite{smid} leads to 
new (infinite dimensional) Hamiltonian loop group
manifolds and $\cC,\wt\cC$ are the corresponding quasi-Hamiltonian spaces.

In more detail recall that the equivalence theorem (Theorem 8.3) of
\cite{AMM} gives a correspondence between Hamiltonian $LK$-manifolds
(with proper moment maps)
and quasi-Hamiltonian $K$-spaces, where $K$ is 
a compact (connected) Lie group and 
$LK=C^\infty(S^1,K)$ is the corresponding loop group.
The main examples of such Hamiltonian $LK$ spaces are moduli spaces of
framed flat connections on 
principal $K$-bundles over compact two-manifolds $\Sigma$
with precisely one boundary component: Given $\Sigma$ and $K$ one
defines a space of connections
$$\A:=\{\al\in \ominfo(\Sigma,\lk)\}$$
on the trivial $C^\infty$ principal $K$-bundle over $\Sigma$ 
(where $\lk=\Lia(K)$) and a
gauge group
$$\cK:= C^\infty(\Sigma,K).$$
This has normal subgroup 
$\cK_\partial:= \{g\in\cK\ \bigl\vert\ g\vert_{\partial\Sigma}=1\}$
consisting of bundle automorphisms equal to the identity on the
boundary circle.
The quotient $\cK/\cK_\partial$ is thus isomorphic to the loop group
$LK$.
Atiyah--Bott \cite{AB83} define the following symplectic structure on
$\A$:
$$\omega_\A(\phi,\psi) = \int_\Sigma(\phi,\psi)$$
where $(,)$ denotes a chosen pairing on $\lk$.
Then, taking  the curvature of the connections in $\A$ gives a moment
map for the action of $\cK_\partial$ (see \cite{Aud95}) and so the
symplectic quotient at the zero value of the moment map is the moduli
space of flat connections with a framing along the boundary circle:
$$\wh\cN := \A_\flt/\cK_\partial.$$
This infinite dimensional symplectic manifold is a Hamiltonian
$LK$-space in the sense of \cite{AMM} (and such spaces 
constitute the main class of examples).
The action of $LK$ is simply the residual action of $\cK$
and the moment map is the restriction of the connections to the
boundary circle:
$$\wh\mu:\wh\cN\to \A_{S^1};\quad\al\mapsto \al\vert_{\partial \Sigma}.$$
(This is really a Hamiltonian $\wh{LK}$-space where $\wh{LK}$ is the
centrally extended loop group and the central circle acts trivially on
$\wh \cN$; the space $\A_{S^1}$ of connections on
the trivial $K$-bundle over the circle is naturally identified with
the level one hyperplane in the dual of the Lie algebra of $\wh{LK}$.
However this complication is incorporated into the definition of
Hamiltonian $LK$-spaces in \cite{AMM, MW-Cob}.)

Now choose a point $p\in\partial \Sigma$ of the boundary circle of $\Sigma$. 
The equivalence theorem of \cite{AMM} implies that the quotient 
$\cN:=\wh \cN/\Omega K$ of $\wh \cN$ by the based loop group 
$\Omega K=\{g\in LK\ \bigl\vert\ g(p)=1\}$ is a (finite dimensional) 
quasi-Hamiltonian $K$-space.
In other words moduli spaces of flat connections on $\Sigma$ with a framing at
one point on the boundary are naturally quasi-Hamiltonian $K$-spaces.

The two-form and moment map on $\cN$ are constructed as follows.
One has a commutative diagram:
\begin{equation}        \label{cd: htoqh}
\begin{array}{ccc}
  \wh\cN & \mapright{\wh\mu} & \A_{S^1} \\
\mapdown{\pi} && \mapdown{h} \\
  \cN & \mapright{\mu} & K.  
\end{array}
\end{equation}
where $\pi$ is the $\Omega K$ quotient and
the maps $\mu$ and $h$ take the holonomy of the connections around
the boundary circle (in a positive sense starting at $p$ with initial
condition $1\in K$). The quasi-Hamiltonian two-form 
$\omega_\cN$ 
on $\cN$ is defined
by\footnote{The signs differ from \cite{AMM} as 1) we give the boundary circle
the induced orientation and 2) an overall sign change has been 
made anyway.}  
$$-\pi^*(\omega_\cN) = \omega_{\wh\cN} - \wh\mu^*(\varpi)$$
where $\omega_{\wh\cN}$ is the symplectic form on $\wh\cN$ and
$\varpi$ is the following two-form on $\A_{S^1}$.
For each point $z\in S^1$ define a map $h_z:\A_{S^1}\to K$
taking a connection $\al$ to its holonomy along the positive arc from 
$p$ to $z$, with initial condition $1\in K$. Thus 
$h_z^*\overline\theta$ is a $z$-dependent 
$\lk$-valued one-form on $\A_{S^1}$ and
$\varpi$ is defined to be 
\begin{equation} \label{eq: varpi def}
\varpi =
\frac{1}{2}\int_{S^1}(h_z^*\overline\theta,dh_z^*\overline\theta)
\end{equation}
where $d$ is the exterior derivative on $S^1$.
It is worth noting that this procedure of subtracting off
$\wh\mu^*(\varpi)$ will amount simply to forgetting part of an integral in
the computation below.

\begin{rmk}
Under this map from surfaces with just one boundary component to
quasi-Hamiltonian $K$-spaces, the quasi-Hamiltonian fusion operation 
corresponds to gluing two surfaces (each with one boundary component)
into two of the holes of a three-holed sphere
(so the resulting surface again has one boundary component)
cf. \cite{AMM, MW-Cob}. 
Also, quasi-Hamiltonian reduction corresponds
to fixing the conjugacy class of monodromy around the boundary
component and forgetting the framing, thereby giving the usual symplectic
moduli space of flat connections.
The upshot is that once we allow fusion, all the symplectic manifolds
that arise as moduli spaces of flat connections on surfaces 
may be constructed from just two
types of quasi-Hamiltonian $K$-spaces: 
conjugacy classes (one for each boundary component)
and the internally fused double ($\cong K\times K$), which corresponds
to the one-holed torus.
\end{rmk}

Now we will apply the above philosophy to the extension of the
Atiyah--Bott symplectic structure to singular connections
($C^\infty$ connections with poles) given in \cite{smid}.
First we point out that the above story may be
complexified; if $\Sigma$ has just one boundary component
and $G$ is a connected complex reductive group (e.g. the complexification
of $K$) then 
the moduli space 
of flat connections on $G$-bundles over $\Sigma$ with framings at one point on
the boundary are complex quasi-Hamiltonian $G$-spaces.
In turn if $\Sigma$ has a complex structure 
such moduli spaces may be
identified with the moduli space of holomorphic connections on
holomorphic $G$-bundles over $\Sigma$ (together with a framing at one point on
the boundary). (Both spaces are isomorphic to the manifold
$\Hom(\pi_1(\Sigma,p),G)$ of fundamental group representations.)

In a similar way the moduli spaces of flat $C^\infty$ singular
connections we will define below
correspond both to moduli spaces of meromorphic connections on holomorphic 
$G$-bundles (cf. \cite{smid} Proposition 4.5) and to spaces of
monodromy/Stokes data (cf. \cite{smid} Proposition 4.8).

Due to fusion it is sufficient to consider only $C^\infty$ singular
connections on a disc having just one pole.
Fix an integer $k\ge 1$ (the pole order) and an {\em irregular type}
$$\wt A^0:= A_0\frac{dz}{z^k}+\cdots+ A_{k-2}\frac{dz}{z^2}\in
\Omega^1[D](\Delta,\g)$$
where $A_i\in \lt, A_0\in\lt_\reg$, $z$ is a coordinate on the closed
unit disc $\Delta$ and $D:=k(0)$ is a divisor on $\Delta$ supported at
the origin.
If $k=1$ we set $\wt A^0=0$.
The spaces of $C^\infty$-singular connections we are interested in
have their full infinite jets of derivatives fixed, except for the
residue term:
$$\wt\A:=\{\al\in\ominfo[D](\Delta,\g)\ \bigl\vert\ 
L_0(\al) = \wt A^0 + \Lambda dz/z \text{ for some } \Lambda\in\lt_k\}$$
where $L_0$ takes the full $C^\infty$ Laurent expansion of $\al$ at
the origin and
$\lt_k=\lt$ if $k\ge 2$ but $\lt_1$ is the affine regular Cartan:
$\lt_1=\{\Lambda\in\lt\ \bigl\vert \ \be(\Lambda)\notin\IZ \text{ for all
roots } \be \}$.
Let
$$\G_T := \{g\in C^\infty(\Delta,G)\ \bigl\vert\ L_0(g)\in T\subset
G\flb z,\overline z\frb\}$$
be the group of bundle automorphisms having Taylor expansion zero at
the origin except for the constant term, which should be in $T$.
Clearly the tangent space to $\wt\A$ at a connection $\al$ is
$$T_\al\wt\A = \{\phi\in\ominfo[D](\Delta,\g)\ \bigl\vert\ 
L_0(\phi)\in\lt \frac{dz}{z}\}.$$
Thus as in \cite{smid} 
we may still use the Atiyah--Bott formula in this singular
situation and define a symplectic structure on $\wt\A$ as
$$\omega_{\wt\A}(\phi,\psi)= \int_\Delta(\phi,\psi).$$
\begin{lem}
The gauge action of the subgroup
$\G_{1,\partial}:=\{g\in\G_T\ \bigl\vert\ g\vert_{\partial\Delta}=1, 
g(0)=1\}$ on $\wt\A$ is Hamiltonian with moment map given by the
curvature.
\end{lem}
\pf
See \cite{smid} Proposition 5.4.
\epf
The symplectic quotient of $\wt\A$ at the zero value of the moment map
is thus
$$\wh\cN := \wt\A_\flt/\G_{1,\partial}$$
which has a residual action of $\G_T/\G_{1,\partial}\cong T\times LG$.
The $T$-action is Hamiltonian with moment map 
$$\al\mapsto -2\pi i\Lambda = - (2\pi i)\Res_0 L_0(\al)$$
as in \cite{smid} Proposition 5.5, and (as above) the $LG$-action
makes $\wh\cN$ into a Hamiltonian $LG$-space (in the sense of
\cite{AMM}) with moment map
$$\wh\mu:\wh\cN \to \A_{S^1}; \quad \al\mapsto\al\vert_{\partial\Delta}.$$
Now fix the point $p=-1\in\partial\Delta$.
Thus (momentarily forgetting the $T$-action)
the quotient $\cN:=\wh\cN/\Omega G$ by the based loop group
should be a quasi-Hamiltonian $G$-space.
First we will use the irregular Riemann-Hilbert correspondence to
identify $\cN$ as a complex manifold.
Let 
$$\G_{1,p}:=\{g\in\G_T\ \bigl\vert\ g(p)=1=g(0)\}$$
so that 
$$\cN = \wh\cN/\Omega G = \wt\A_\flt/\G_{1,p}$$
which has a residual action of $\G_T/\G_{1,p}\cong G\times T$.
\begin{thm}[\cite{smid,bafi}] \label{thm: what c is}
The quotient $\wt\A_\flt/\G_{1,p}$ is isomorphic to $\wt \cC$ as a
$G\times T$-space.
\end{thm}
\pf
As in \cite{smid} Proposition 4.5, Corollary 4.6 this quotient may be
shown to be canonically isomorphic to the set of isomorphism classes of 
$4$-tuples $(P,A,g_0,g_p)$ where $P\to \Delta$ is a holomorphic
principal $G$-bundle, $A$ is a meromorphic connection on $P$ with
irregular type $\wt A^0$ and compatible framing $g_0$ at the origin
and $g_p$ is an arbitrary framing of $P$ at $p$.
Then by the irregular Riemann-Hilbert correspondence of \cite{bafi} Section 2 
the
moduli space of such triples $(P,A,g_0)$ 
is analytically isomorphic to the space 
$(U_+\times U_-)^{k-1}\times\lt$ of Stokes multipliers and exponents
of formal monodromy ($\Lambda$'s). The inclusion of the framing $g_p$
in the moduli problem simply adds a factor of $G$ so the result
follows.
The formula for the $G$-action is immediate and for the $T$-action see
\cite{smid} Corollary 3.5.
\epf

\begin{rmk}
The  monodromy map
$\wt \nu:\wt\A_\flt\twoheadrightarrow\wt\cC$, whose fibres are
precisely the $\G_{1,p}$ orbits will be described directly in the proof of
the following theorem.
\end{rmk}

Now if $\omega_{\wh\cN}$ is the symplectic structure on $\wh\cN$ and 
$\varpi$ is the complex analogue of the two-form 
\eqref{eq: varpi def} on $\A_{S^1}$ (defined exactly the same way) then, 
by the general theory described above, we expect the the two-form 
$-\omega_{\wh\cN}+\wh\mu^*(\varpi)$ on $\wh\cN$ to be the pullback of some
quasi-Hamiltonian two-form on $\wt\cC$ along the 
map $\pi:\wh\cN\to\cN\cong\wt\cC$. Indeed we have the following
theorem.

\begin{thm} \label{thm: main derivation}
Let $\omega$ be the two-form on $\wt \cC$ defined in 
\eqref{eq: first omega exprn}. Then we have
$$-\pi^*(\omega)=\omega_{\wh\cN}-\wh\mu^*(\varpi).$$
\end{thm}

\pf 
Since $\wh\cN$ is the symplectic quotient of $\wt\A$ this is
equivalent to proving 
$\iota^*\omega_{\wt\A}-\pr^*\wh\mu^*\varpi=-\wt\nu^*\omega$ where
$\iota:\wt\A_\flt\to\wt\A$ is the inclusion, $\pr:\wt\A_\flt\to\wh\cN$
is the projection and $\wt\nu:\wt\A_\flt\to\wt\cC$.
To this end suppose we have a two-parameter family
$\al(s,t)\in\wt\A_\flt$ of flat singular connections depending
holomorphically on $s,t$.
We will evaluate the two-form  $\iota^*\omega_{\wt\A}-\pr^*\wh\mu^*\varpi$
on the pair $\al',\dot\al\in\ominfo[D](\Delta,\g)$ of
tangent vectors to $\wt\A_\flt$ at $\al=\al(0,0)$, where
$\al'=\frac{d}{ds}\al\bigl\vert_{s=t=0}, 
 \dot\al=\frac{d}{dt}\al\bigl\vert_{s=t=0}$.
If $X=\wt\nu_*(\al'), Y=\wt\nu_*(\dot\al)\in T_{\wt\nu(\al)}\wt\cC$ then
we should obtain $-\omega(X,Y)$ where by definition
$$2\omega(X,Y) = 
(\overline\cD',\dot{\overline\cE}) -
(\dot{\overline\cD},\overline\cE') + 
\sum_{j=1}^{k-1} 
(\cD'_j,\dot\cD_{j-1})- (\dot\cD_j,\cD'_{j-1})
- (\cE'_j,\dot\cE_{j-1})+(\cE'_j,\dot\cE_{j-1})
$$
with $\cD'_j = \langle\cD_j,X\rangle$ etc.

Let $\Delta_r$ denote the slit annulus obtained by cutting $\Delta$ along the
ray from $0$ to $p=-1$ and removing the open disc of radius $r$ centred on the
origin.
Denote by $\csan$ the closure of $\Delta_r$ in the universal cover of the
punctured disc $\Delta\setminus \{0\}$.
Thus $\csan$ has two straight edges $l_+,l_-$
lying over the interval $[-1,-r]$ and has 
interior isomorphic to the interior of 
$\Delta_r\subset \Delta$.
In particular $\csan$ is contractible.
We identify the lower lip $l_-$ with  the interval $[-1,-r]\subset \Delta$,
so that one arrives at the upper lip $l_+$ by turning a full turn in a
positive sense from $l_-$. 
For each $s,t$ let
$$\chi(s,t):\csan\to G$$
be the fundamental solution of the connection $\al(s,t)$ taking
the value $1\in G$ at $p\in l_-$. (In other words $\chi(s,t)$ is
the 
map solving the differential equation $\al(s,t)=\chi^*(\overline\theta)$.)
Then for each $z\in\csan$ let
$\chi'(z):=\frac{d}{dt}\chi(s,t,z)\bigl\vert_{s=t=0}\in T_{\chi(z)}G$
and so $$\chi^{-1}\chi' :=l_{\chi^{-1}}\chi'$$ is a $\g$-valued function
on $\csan$, where $l_{\chi^{-1}(z)}:T_{\chi(z)}G\to \g$ 
denotes the derivative of left
multiplication by $\chi^{-1}(z)$ in the group $G$.
Now define a one-form $\varphi$ on $\csan$ by
$$\varphi:= 
\frac{1}{2}(\varphi_1-\varphi_2), \qquad
\varphi_1:=
(\chi^{-1}\chi',d(\chi^{-1}\dot\chi)),\quad
\varphi_2:= (d(\chi^{-1}\chi'),\chi^{-1}\dot\chi)$$
where $d$ is the exterior derivative on $\csan$.
Thus $d\varphi=(\al',\dot\al)$ as two-forms on $\csan$ (since 
e.g. $\al'=\chi d(\chi^{-1}\chi')\chi^{-1}$).
In turn since $(\al',\dot\al)$ is a smooth two-form on $\Delta$ we
have
$$\omega_{\wt\A}(\al',\dot\al)=\int_\Delta(\al',\dot\al)=
\lim_{r\to 0}\int_{\csan}d\varphi = 
\lim_{r\to 0}\int_{\partial\csan} \varphi.$$
This integral will be evaluated along each arc of the boundary of 
$\csan$, neglecting any terms that vanish in the limit.
A similar calculation appears in 
\cite{Nick-sat}.

First around the outer boundary of $\csan$ (the circle of radius one) 
we recognise that the integral of $\varphi$ is
precisely $(\pr^*\wh\mu^*\varpi)(\al',\dot\al)$ (since on this circle $\chi$
restricts to the map $h_z$ used to define $\varpi$),
which is the term to be subtracting off.

For the other arcs we first need to describe directly the map
$\wt\nu:\wt\A_\flt\to\wt\cC$ 
associating monodromy data $(C,\bfd,\bfe)$ to a flat singular connection $\al$.
The key point is that any $\al\in\wt\A_\flt$ has canonical
fundamental solutions
$$\Phi_i:\Sect_i\to G$$
on certain distinguished sectors $\Sect_i$ defined as follows
(\cite{smid} Lemma 4.7, \cite{bafi} Section 2).
The leading coefficient $A_0\in\lt_\reg$ of the chosen irregular type
$\wt A^0$ determines the {\em anti-Stokes directions} $\IA$ at
$0\in\Delta$ defined as
$$z\in\Delta\setminus\{0\}\text{ 
lies on an anti-Stokes direction }\Longleftrightarrow
\frac{\be(A_0)}{z^{k-1}}\in\IR\text{ 
for some root }\be\in\R.$$
This determines a finite set $\IA$ of directions which is clearly 
invariant under rotation by $\pi/(k-1)$ and so the number
$l:=\#\IA/(2k-2)$ is an integer.
The sectors $\Sect_i$ are just the sectors bounded by consecutive
anti-Stokes directions.
Without loss of generality we will assume the positive real axis $\IR_+$
is not an anti-Stokes direction and label these sectors in a positive
sense and such that $\IR_+\subset \Sect_0$.
In turn the anti-Stokes directions $a_i\in\IA$ 
are labeled (modulo $\#\IA$) such
that $\Sect_i=\Sect(a_i,a_{i+1})$.
By \cite{bafi} Lemma 2.4 we know that the set of roots
$$\R_+:= \{\be\in\R\ \bigl\vert\ \frac{\be(A_0)}{z^{k-1}}\in\IR_+
\text{ for $z$ on one of the directions $a_1,\ldots,a_l$ }\}$$
`supporting' one of the first $l$ anti-Stokes directions, is a set of
positive roots, and we define $B_+$ to be the corresponding Borel subgroup
containing $T$. Now to define $\Phi_i$ we recall that the Laurent
expansion of $\al$ is
$$L_0(\al)=dQ+\Lambda\frac{dz}{z}=:A^0$$
for some $\Lambda\in\lt$ where
$Q:=\sum_{j=1}^{k-1}\frac{z^{j-k}}{j-k}A^0_{j-1}$ (so $dQ=\wt A^0$).
In particular the $(0,1)$ part of $\al$ is nonsingular across the
origin and so we may solve the $\overline\partial$-problem
$$(\overline\partial g)g^{-1} = \al^{0,1}$$
for a smooth map $g:U\to G$ defined in some neighbourhood
$U\subset\Delta$ of the origin.
Given such $g$ one observes (\cite{smid} Lemma 4.3) that the Taylor
expansion $\wh F=L_0(g^{-1})$ is in $G\flb z\frb$ (has no $\overline
z$ terms) and that $A:=\wh F[A^0]=g^{-1}[\al]$ is the germ of a
(convergent) meromorphic connection.
In turn this implies (\cite{bafi} Theorem 2.5) that there is a unique
holomorphic map
$$\Sigma_i(\wh F):\Sect_i\to G$$
on each sector such that
$\Sigma_i(\wh F)[A^0]=A$ and that the analytic continuation of 
$\Sigma_i(\wh F)$ to the {\em supersector}
$$\Ssect_i:=\Sect\left(a_i-\frac{\pi}{2k-2},a_{i+1}+
\frac{\pi}{2k-2}\right)$$
is asymptotic to $\wh F$ at $0$ in $\Ssect_i$.
Now we are led to the
following definition because 
$z^\Lambda e^Q$ is a fundamental solution of the connection $A^0$,
$\Sigma_i(\wh F)$ is an isomorphism between $A^0$ and $A$, and $g$ is
an isomorphism between $A$ and $\al$.
(Here $z^\Lambda$ is defined on $\Sect_0$ using the branch of
$\log(z)$ that is real on $\IR_+$ and by convention we extend this to
the other sectors in a {\em negative} sense.)
\begin{defn}
The {\em canonical fundamental solution} of $\al\in\wt\A_\flt$ on
$\Sect_i$ is the map
$$\Phi_i:= g\Sigma_i(L_0g^{-1})z^\Lambda e^Q:\Sect_i\to G$$
for {\em any} solution $g$ of $(\overline\partial g)g^{-1} = \al^{0,1}$.
\end{defn}
The Stokes multipliers $S_i$ of $\al$ can now be defined (as in \cite{bafi}
Definition 2.6) as the elements of $G$ relating the fundamental solutions 
$\Phi_{il}$ and $\Phi_{(i+1)l}$.
However to define directly the elements $d_i,e_i$ we first define new
fundamental solutions $\Psi_i, \Theta_i$ as follows:
$$\Psi_i:= \Phi_{il}\epsilon^{2k-2-i}:\Sect_{il}\to G
\quad(i=1,\ldots,2k-2), \qquad\Theta_i=\Psi_{2k-2-i}.$$
where $\epsilon:=e^{\frac{\pi i\Lambda}{k-1}}$.
The indices of $\Psi_i, \Theta_i$ are taken modulo $2k-2$ so
$\Psi_0=\Theta_0=\Phi_0$ on $\Sect_0$.
For $i=0,\ldots,k-2$ the sector on which $\Psi_i$ or $\Theta_i$ is
defined intersects the slit annulus 
$\Delta_r$ in a contractible set and so we may
extend $\Psi_i,\Theta_i$ uniquely (as fundamental solutions of  $\al$)
to maps  from $\csan$ to $G$.
Now the intersection of 
$\Sect_{(k-1)l}$ (the sector
containing $\IR_-$) and $\Delta_r$ has two components, and we extend 
$\Psi_{k-1}$ from the upper component of this intersection  onto
$\csan$ and we extend $\Theta_{k-1}$ from the lower component.
Thus we have $2k$ generally distinct 
fundamental solutions of $\al$ on $\csan$:
$$\chi,\Phi_0=\Psi_0=\Theta_0,
\Psi_1,\ldots,\Psi_{k-1},\Theta_1,\ldots,\Theta_{k-1}.$$
The monodromy data $C,d_i,e_i$ 
is defined to be the set of ($z$-independent) group
elements relating them, as follows:
$$\Phi_0C=\chi,\quad
\Psi_ie_i=\Psi_{i-1},\quad\Theta_id_i=\Theta_{i-1}\quad(i=1,\ldots,k-1).$$
If $d_i,e_i$ are defined in this way it follows from \cite{bafi} Lemma
2.7 that
$d_{\text{even}} , e_{\text{odd}}\in B_+$, 
$d_{\text{odd}},e_{\text{even}}\in B_-$ and
$\delta(d_j)^{-1} = \epsilon = \delta(e_j)$, so we have indeed
associated a point of $\wt\cC$ to $\al$.
Note also that the maps $D_i,E_i:\wt\cC\to G$ arise as
$$\Psi_iE_i=\chi,\quad\Theta_iD_i=\chi\quad(i=0,\ldots,k-1).$$
It follows that $\chi$ has holonomy $D^{-1}E$ since 
$\chi\vert_{l_+}=\Psi_{k-1}\vert_{l_+}E=\Theta_{k-1}\vert_{l_-}E=
\chi\vert_{l_-}D^{-1}E$,
and so this is the 
quasi-Hamiltonian monodromy map.

Now we return to the boundary integral.
Choose a point $q_i$ of distance $r$ from the origin and in the intersection 
$\Ssect_{il}\cap\Ssect_{(i-1)l}$ of two of the supersectors, for 
$i=1,\ldots,k-1$. Thus we know that both $\Phi_{il}(q_i)$ and
$\Phi_{(i-1)l}(q_i)$ are asymptotic to $z^\Lambda e^Q$ at $0$ as $r\to 0$, and
in turn we know the asymptotics of $\Psi_i(q_i)$ and
$\Psi_{i-1}(q_i)$ at $0$.
Similarly choose $p_i\in \Ssect_{-il}\cap\Ssect_{-(i-1)l}$ of modulus $r$ so
that we know the asymptotics of both $\Theta_i(p_i)$ and
$\Theta_{i-1}(p_i)$ at $0$ as $r\to 0$.
Let $p_k=-r\in l_-$ and let $q_k$ be the point of the upper lip $l_+$ lying
over $-r$. 
Thus we 
may divide the inner boundary circle of $\csan$ into 
$2k-1$ arcs by breaking it at the points $p_i,q_i$.
Now since $\chi=\Psi_iE_i$ and $E_i$ is $z$-independent we find
\begin{equation} \label{eq: phi1 expn}
\varphi_1=(\Psi_i^{-1}\Psi_i',d(\Psi_i^{-1}\dot\Psi_i)) + 
d(\overline\cE'_i,\Psi_i^{-1}\dot\Psi_i)
\end{equation}
where $\overline\cE'_i=\langle E_i^*\overline\theta,X\rangle$, and similarly
for $\varphi_2$ (swapping the dot and the prime).
\begin{lem}
The first term in \eqref{eq: phi1 expn} 
may be neglected in the integral from $q_{i+1}$ to $q_i$.
\end{lem}
\pf
The first term of \eqref{eq: phi1 expn} and the corresponding term of
$\varphi_2$ contribute  
\begin{equation} \label{eq: phi int}
\frac{1}{2}\int_{q_{i+1}}^{q_i}(\Psi_i^{-1}\Psi_i',d(\Psi_i^{-1}\dot\Psi_i))-
(\Psi_i^{-1}\dot\Psi_i,d(\Psi_i^{-1}\Psi'_i))
\end{equation}
to the integral of $\varphi$.
However $\Psi_i\simeq z^\Lambda e^Q\epsilon^{2k-2-i}$ at $0$ uniformly in
$\Ssect_{il}$ (which contains the integration path). 
Substituting in this approximation gives that the integrand 
in \eqref{eq: phi int} is zero. This implies that
in the limit $r\to 0$ the integral 
\eqref{eq: phi int} really is zero.
\epf

Thus modulo negligible terms
$$\int_{q_{i+1}}^{q_i}\varphi_1=
(\overline\cE'_i,\Psi_i^{-1}\dot\Psi_i)\bigl\vert_{q_{i+1}}^{q_i}.$$
If we sum this integral for $i=1,\ldots,k-1$ then the contribution at $q_i$ is
$$(\overline\cE'_i,\Psi_i^{-1}\dot\Psi_i)(q_i)-
(\overline\cE'_{i-1},\Psi_{i-1}^{-1}\dot\Psi_{i-1})(q_i)$$
provided  $i\ne 1, k$.
Now using $\Psi_{i-1}=\Psi_ie_i$ to remove $\Psi_{i-1}$ this becomes
$$(\overline\cE'_i-e_i\overline\cE'_{i-1}e_i^{-1} , \Psi_i^{-1}\dot\Psi_i)-
(\overline\cE_{i-1}',\dot\varepsilon_i)$$
where $\dot\varepsilon_i=\langle e_i^*\theta,Y\rangle$.
In turn using $E_i=e_iE_{i-1}$ this becomes 
\begin{equation} \label{eq fianl phi1 exprn}
(\varepsilon'_i , \Psi_i^{-1}\dot\Psi_i)-
(\cE_{i-1}',\dot\cE_i).
\end{equation}
If we also repeat the above for $\varphi_2$ we get the same but with 
the dots and primes swapped.
Now 
since $\Psi_i\simeq z^\Lambda e^Q\epsilon^{2k-2-i}$ and the $T$ component of
$e_i$ is $\epsilon$ we deduce 
$(\varepsilon'_i , \Psi_i^{-1}\dot\Psi_i)-
 (\dot \varepsilon_i , \Psi_i^{-1}\Psi'_i) \to 0$ as $r\to 0$.
Thus the contribution at $q_i$ ($i\ne 1,k$) 
to the integral of $\varphi$ from $q_k$ to $q_1$ 
is
$$-\frac{1}{2}((\cE_{i-1}',\dot\cE_i)-(\dot\cE_{i-1},\cE'_i)) = 
\frac{1}{2}(\cE_{i},\cE_{i-1})(X,Y)$$
which is a term appearing in $-\omega(X,Y)$.
Writing $p_0:=q_1$ and performing the same manipulations for the $\Theta_i$'s,
integrating $\varphi$ from $p_0$ to $p_k$ yields a contribution of 
$$\frac{1}{2}((\cD_{i-1}',\dot\cD_i)-(\dot\cD_{i-1},\cD'_i)) = 
-\frac{1}{2}(\cD_{i},\cD_{i-1})(X,Y)$$
at $p_i$, provided $i\ne 0,k$.
The two left-over contributions at $q_1=p_0$ combine to give the term 
$\frac{1}{2}(\cE_{1},\cE_{0})(X,Y)$. (Thus all terms of $-\omega$ except
$-\frac{1}{2}(\overline\cD,\overline\cE)(X,Y)$ have been obtained so far.) The
left-over contributions at $q_k$
and $p_k$ are:
$$\frac{1}{2}((\overline\cD',\Theta_{k-1}^{-1}\dot\Theta_{k-1})-
(\dot{\overline\cD},\Theta_{k-1}^{-1}\Theta'_{k-1}))(p_k)
-\frac{1}{2}((\overline\cE',\Psi_{k-1}^{-1}\dot\Psi_{k-1})-
(\dot{\overline\cE},\Psi_{k-1}^{-1}\Psi'_{k-1}))(q_k).$$

Now consider the two straight edges $l_\pm$ of $\csan$.
Recall that $\Theta_{k-1}\vert_{l_-}=\Psi_{k-1}\vert_{l_+}$ so that 
from \eqref{eq: phi1 expn}
$$\int_{l_++l_-}\varphi_1 = 
\int_{p_k}^p d(\overline\cD'-\overline\cE',\Theta_{k-1}^{-1}\dot\Theta_{k-1})=
(\overline\cD'-\overline\cE',
\Theta_{k-1}^{-1}\dot\Theta_{k-1})\bigl\vert_{p_k}^p
$$
and similarly for $\varphi_2$. 
Observe that 
the  contribution at $p_k$ to the integral of $\varphi$ along $l_\pm$ 
cancels precisely with the left-over terms at $p_k,q_k$ displayed above.
Finally since $\Theta_{k-1}(p)=\chi(p)D^{-1}=D^{-1}$ the contribution at $p$ is
$$ -\frac{1}{2}(
(\overline\cD'-\overline\cE',\dot{\overline \cD})-
(\dot{\overline\cD}-\dot{\overline\cE},\overline \cD'))
=-\frac{1}{2}(\overline\cD,\overline\cE)(X,Y).$$

\epf
\end{section}

\section{Additive Analogues}  \label{sn: additive}
Here we recall (from \cite{smid} Section 2) the symplectic manifolds 
$O,\wt O$ which are the additive analogues of the quasi-Hamiltonian
spaces $\cC,\wt \cC$. 

Fix an integer $k\ge 2$. 
Let $G_{k}:= G\bigl(\IC[z]/z^{k}\bigr)$
be the group of $(k-1)$-jets of bundle automorphisms, and let 
$\g_k=\Lia(G_k)$ be its Lie algebra, which contains elements of the
form 
$X=X_0+X_1z+\cdots+X_{k-1}z^{k-1}$ with $X_i\in\g$.
Let $B_k$ be the subgroup of $G_k$ of elements having constant term
$1$. The group $G_k$ is the
semi-direct product  $G\ltimes B_k$ (where 
$G$ acts on $B_k$ by conjugation).
Correspondingly the Lie algebra of $G_k$ decomposes as a vector space
direct sum and dualising we have:
$\gks=\lbks\oplus \gs$.
Elements of $\gks$ will be written as 
\begin{equation} 	\label{eqn: pp map}
A=A_{0}\frac{dz}{z^{k}}+\cdots+ A_{k-1}\frac{dz}{z}
\end{equation}
via the pairing with $\g_k$ given by
$\langle A,X\rangle := \res_0(A,X)=
\sum_{i+j=k-1} (A_i,X_j)$.
In this way 
$\lbks$ is identified with the set of $A$ 
having zero residue and $\g^*$ with those having only a residue term 
(zero irregular part).
Let $\pir:\gks\to\gs$ and $\pii:\gks\to\lbks$ denote the corresponding
projections.

Now 
choose an element 
$\wt A^0=
A^0_0dz/z^k+\cdots+A^0_{k-2}dz/z^{2}$
of $\lbk^*$ with $A^0_i\in\lt$ and with regular leading coefficient 
$A^0_0\in\lt_\reg$.
Let $O_B\subset \lbks$ denote the $B_k$ coadjoint orbit containing
$\wt A^0$.

\begin{defn}	\label{dfn: ext orb}
The {\em extended orbit} $\wt O\subset G \times \gks$ associated to $O_B$ is:
$$\wt O := 
\left\{
(g_0,A)\in G \times \gks \ \bigl\vert
\ \pii(g_0 A g_0^{-1})\in O_B \right\} $$ 
where $\pii:\g_k^*\to\lbks$ is the natural projection removing the residue.
\end{defn}
If $(g_0,A)\in\wt O$ then $A$ 
will correspond to the principal part of a generic
meromorphic connection and $g_0$ to a compatible framing.

In the simple pole case $k=1$ we define 
$$
\wt O := \left\{ (g_0,A)\in G\times \gs
 \ \bigl\vert \ g_0Ag_0^{-1}\in \lt_1 \right\}\subset G\times \gs
$$
where $\lt_1\subset \lt^*\cong\lt$ is the complement of the affine
root hyperplanes.
If we identify $G\times \gs$ with $T^*G$ then $\wt O$ is in fact a
symplectic submanifold (see \cite{GS} Theorem 26.7).

The basic properties of these extended orbits may be summarised as follows. 
Given $(g_0,A)\in\wt O$ then by hypothesis there is some 
$g\in G_k$ such that  $gAg^{-1} = \wt A^0 + R dz/z$ for some 
$R\in \g$ and we define a map $\Lambda=\delta(R):\wt O\to
\lt\cong\lt^*$ 
by taking the $\lt$ component of $R$ (which is independent of $g$). 

\begin{prop}[\cite{smid}] 

1). 
The extended orbit $\wt O$ is canonically isomorphic to the symplectic
quotient $(T^*G_k\times O_B)\spq B_k$.

2). (Decoupling).
The map $\wt O \to (T^*G) \times O_B; 
(g_0,A)\mapsto(g_0,\pir(A), \pii(g_0Ag_0^{-1}) )$
is a symplectic isomorphism
where $T^*G\cong G\times\gs$ via the left trivialisation.

3). 
The map $-\Lambda$ is a moment map for 
the  free action of $T$ on $\wt O$ defined by
$t(g_0,A)=(tg_0,A)$ where $t\in T$.

4). 
The symplectic quotient by $T$ at the value $-\Lambda$ of the moment
map is the $G_k$ coadjoint orbit $O$ through the element 
$\wt A^0+\Lambda dz/z$ of $\gks$.

5). 
The free $G$-action $h(g_0,A):=(g_0h^{-1},hAh^{-1})$ on $\wt O$ 
is Hamiltonian with moment map
$\mu_{G}: \wt O \to \gs; (g_0,A)\mapsto \pir(A).$
\end{prop}

In particular $\wt O$ is a Hamiltonian $G\times T$-manifold with $T$
reductions equal to $G_k$-coadjoint orbits $O$; these
properties are viewed as natural analogues of those of $\wt \cC$ (and
they do indeed match up under the Riemann-Hilbert correspondence). 
Note that the coadjoint orbit $O_B$ is a point if $k=2$ so that 
part 2) says $\wt O\cong T^*G$, the additive analogue of the fact that
$\wt\cC\cong G\times G^*$ in this case.

Proposition 2.1 of \cite{smid}
explains how the symplectic manifolds 
$$(\wt O_1\times\cdots\times \wt O_m)\spq G$$
for extended orbits $\wt O_i$,
are isomorphic to moduli spaces of (compatibly framed)
meromorphic connections on trivial $G$-bundles over $\IP^1$ (with fixed
irregular types).

\appendix
\begin{section}{Kernel Calculation}

We will establish the minimal degeneracy condition (QH3) for the two-form
$\omega$ on $\wt \cC$. 

\ 

\pfms {\em(of (QH3)).}\ \ 
\begin{lem}     \label{lem: long 2form}
The two-form $2\omega$ on $\wt \cC$ is also given by the formula
$$\Bigl(\overline\gamma,(11)\overline\gamma(11)^{-1}\Bigr)+
\sum_{i=1}^{k-1}
\Bigl(\overline\gamma, \we(1i)\varepsilon_i(1i)^{-1}+ 
\{i 1\}^{-1}\varepsilon_i\we\{i 1\}-
(1i)^{-1}\delta_i(1i)- [i 1]^{-1}\delta_i\we[i 1]\Bigr)$$
$$\
+\sum_{1\le i,j \le k-1} 
\Bigl(\delta_i,(ij)\varepsilon_j(ij)^{-1}\Bigr) +
\sum_{1\le j < i \le k-1}
\Bigl(\delta_i, [i j]\delta_j[i j]^{-1}\Bigr) - 
\Bigl(\varepsilon_i,\{i j\}\varepsilon_j\{i j\}^{-1}\Bigr)
$$
where $\delta_i = d_i^*(\theta), \varepsilon_i = e_i^*(\theta),
\overline\gamma= C^*(\overline\theta)$,
$(ij):= 
d_i^{-1}d_{i+1}^{-1}\cdots d_{k-1}^{-1}e_{k-1}\cdots e_{j+1}e_j$,
$[ij]:= d_{i-1}\cdots d_{j}$ and 
$\{ij\}:=e_{i-1}\cdots e_{j}$. 
\end{lem}
\pf
This is a straight-forward 
direct calculation expanding each term in \eqref{eq: 2form}.
\epf

Now suppose we choose a pair of tangent vectors $X,Y$ to $\wt \cC$ at
some point $p$, such that $X$ is in the kernel of $\omega_p$ and $Y$
is arbitrary.
We will use dots/primes to denote derivatives along Y/X respectively,
so e.g.
$\dot\delta_i = \langle Y, \delta_i \rangle \in \g$ and
$\varepsilon'_j = \langle X, \varepsilon_j\rangle \in \g$
(and in any representation of $G$ we have
$\dot\delta_i = d_i^{-1}\dot d_i$ etc).
Our aim is to prove $\delta'_i=\varepsilon'_i=0$ for all $i$ (so $X$
is tangent to the $G$ action) and that 
$\Ad_{\mu(p)}(\gamma') = -\gamma'$, which is the required degeneracy condition.
The equation expressing the fact that 
$X$ is in the kernel of $\omega_p$ is equivalent to
\begin{equation} \label{eqn: main ker eqn}
2\omega(Y,X) = 
(\dot\gamma,\Gamma) + \sum_{i=1}^{k-1}
(\dot\delta_i,\Delta_i) + (\dot\varepsilon_i,\xi_i) = 0
\end{equation}
for all $Y$ 
where $\Gamma, \Delta_i, \xi_i\in\g$ are the corresponding coefficients
involving just $X$ derivatives; explicitly from 
Lemma \ref{lem: long 2form}: 
$$\Delta_i= (i1)\overline\gamma'(i1)^{-1}+[i1]\overline\gamma'[i1]^{-1}+
\sum_{j=1}^{k-1}
(ij)\varepsilon'_j(ij)^{-1}+ 
\sum_{j<i}
[ij]\delta'_j[ij]^{-1}-
\sum_{j>i} 
[ji]^{-1}\delta'_j[ji],$$
$$\xi_i= 
-(1i)^{-1}\overline\gamma'(1i)-\{i1\}\overline\gamma'\{i1\}^{-1}
-\sum_{j=1}^{k-1}
(ji)^{-1}\delta'_j(ji)- 
\sum_{j<i}
\{ij\}\varepsilon'_j\{ij\}^{-1}+
\sum_{j>i} 
\{ji\}^{-1}\varepsilon'_j\{ji\},$$
$$\Gamma= (11)\overline\gamma'(11)^{-1}-(11)^{-1}\overline\gamma'(11)+
\sum_{j=1}^{k-1}
(1j)\varepsilon'_j(1j)^{-1}+ 
\{j1\}^{-1}\varepsilon'_j\{j1\}-
(j1)^{-1}\delta'_j(j1)-
[j1]^{-1}\delta'_j[j1].$$
Since $Y$ is arbitrary \eqref{eqn: main ker eqn} implies
$\Gamma=0$
and (since $(,)$ pairs opposite Borels)
that the piece of $\Delta_i$ in the unipotent subalgebra opposite the
Borel containing $\dot\delta_i$ is zero (and similarly 
the piece of $\xi_i$ in the unipotent subalgebra opposite the
Borel containing $\dot\varepsilon_i$ is zero).
The only other information about $X$ in \eqref{eqn: main ker eqn}
concerns the $\lt$ components as follows.
Since $Y$ is tangent to $\wt \cC$ we have 
$\delta(\dot\varepsilon_j)=-\delta(\dot\delta_j)=\pi i \dot\Lambda\in\lt$
for all $j$, where $\delta:\g\to\lt$ is the projection along the root spaces.
Thus, as $\dot \Lambda$ is arbitrary, \eqref{eqn: main ker eqn} implies
\begin{equation} \label{eqn: t cmpt}
\sum_{i=1}^{k-1} \delta(\Delta_i) = 
\sum_{i=1}^{k-1}\delta(\xi_i)
\end{equation}
where $\delta:\g\to\lt$.
Now we will proceed to deduce the required result.
From the formula for $\Delta_i$ it follows that
$d_i\Delta_id_i^{-1}-\Delta_{i+1}=-\overline\delta'_i-\delta'_{i+1}.$
Thus if we define 
$T_i := d_i\Delta_id_i^{-1}+\delta'_{i+1}=\Delta_{i+1}-\overline\delta'_i$
(for $i=1,\ldots,k-2$)
then the restrictions on the unipotent pieces of
$\Delta_i,\Delta_{i+1}$ imply
$T_i= \delta'_{i+1}-\overline\delta'_i + H_i$
for some $H_i\in\lt$, and so in turn
$$\Delta_i = -\delta'_i + d_i^{-1}H_id_i,\qquad
\Delta_{i+1} = \delta'_{i+1} + H_i.$$
Thus $\Delta_i = \delta'_{i} + H_{i-1}= -\delta'_i + d_i^{-1}H_id_i$ 
so that
$ 2\delta'_i = d^{-1}_iH_id_i - H_{i-1}$
for $i=2,\ldots,k-2$.
Taking the $\lt$ component of this implies
$H_{i-1} = H + H_i$
where $H:= (2\pi i)\Lambda'$ 
(so $H=-2\delta(\delta'_j)=2\delta(\varepsilon'_j)$ for all $j$).
If we {\em define} $H_{k-1}:= \delta(\Delta_{k-1})-H/2$ then 
$\delta(\Delta_i)= H/2 + H_i$ for all $i$, and
since $H_i= (k-1-i)H+H_{k-1}$
this implies 
$$\sum_{i=1}^{k-1} \delta(\Delta_i) = (k-1)H_{k-1}+(k-1)^2H/2.$$ 
A similar exercise in terms of the $\varepsilon_i$ and  $\xi_i$ yields
analogous formulae with some sign changes:
$Y_i := e_i\xi_ie_i^{-1}-\varepsilon'_{i+1}=\xi_{i+1}+
\overline\varepsilon'_i$
(for $i=1,\ldots,k-2$), so that
$Y_i= \overline\varepsilon'_i -\varepsilon'_{i+1}+ K_i$
for some $K_i\in\lt$, and in turn
$ 2\varepsilon'_i = -e^{-1}_iK_ie_i + K_{i-1}$
for $i=2,\ldots,k-2$.
Similarly this implies
$K_{i-1} = H + K_i$ and then
$\delta(\xi_i)= H/2 + K_i$ for all $i$ so that 
$\sum_{i=1}^{k-1} \delta(\xi_i) = (k-1)K_{k-1}+(k-1)^2H/2.$
Thus equation \eqref{eqn: t cmpt} is equivalent to
$H_{k-1} = K_{k-1}$.

Now we will reconsider the equations 
$\Delta_i + \delta'_i  = d_i^{-1}H_id_i$ and
$\xi_i-\varepsilon'_i = e_i^{-1}K_ie_i$. 
Using these and the initial formulae for $\Delta_i, \xi_i$ one finds 
$[i1]^{-1}(\Delta_i + \delta'_i)[i1]
+(1i)(\xi_i-\varepsilon'_i)(1i)^{-1}$
is equal to both sides of
\begin{equation} \label{eqn: twisted sum}
2\sum_{j>i}\bigl(
(1j)\varepsilon'_j(1j)^{-1} -  [j1]^{-1}\delta'_j[j1]\bigr)=
[i1]^{-1}(d_i^{-1}H_id_i)[i1] + (1i)(e_i^{-1}K_ie_i)(1i)^{-1}
\end{equation}
Conjugating by $(11)^{-1}$ this is equivalent to
\begin{equation} \label{eqn: twisted sum2}
2\sum_{j>i}\bigl(
\{j1\}^{-1}\varepsilon'_j\{j1\} -  (j1)^{-1}\delta'_j(j1)\bigr) =
(i+11)^{-1}H_i(i+11) + \{i+11\}^{-1}K_i\{i+11\}.
\end{equation}
Putting $i=k-2$ in \eqref{eqn: twisted sum} (so the sum has just one
term) we find
$$2\overline\varepsilon'_{k-1}-2\overline\delta'_{k-1} = 
d_{k-1}H_{k-2}d^{-1}_{k-1} +  e_{k-1}K_{k-2}e_{k-1}^{-1}.$$
Firstly the $\lt$ component of this says
$2H = H_{k-2} + K_{k-2}$ but 
$H_{k-2} = H+H_{k-1}= H+K_{k-1} = K_{k-2}$ and thus we deduce
$H_{k-1} = K_{k-1}=0$ (so that now $H_i=K_i=(k-1-i)H$ for all $i$).
Secondly, rewriting gives 
$$2\overline\varepsilon'_{k-1}-e_{k-1}K_{k-2}e_{k-1}^{-1}=
2\overline\delta'_{k-1} +d_{k-1}H_{k-2}d^{-1}_{k-1}$$
the two sides of which live in opposite Borel subalgebras and 
have zero $\lt$ component, and so are both zero, i.e.
$\varepsilon'_{k-1}=H/2=-\delta'_{k-1}.$

Similarly, considering the difference $[i1]^{-1}(\Delta_i + \delta'_i)[i1] -
   (1i)(\xi_i-\varepsilon'_i)(1i)^{-1}$
instead and setting $i=1$, one obtains 
\begin{equation} \label{eqn: twisted diff result}
2(11)\varepsilon'_1(11)^{-1}+
2(11)\overline\gamma'(11)^{-1} =
-2\delta'_1 -2\overline\gamma'+
(k-2)(d_1^{-1}Hd_1 - (12)H(12)^{-1}).
\end{equation}
Conjugating by $(11)^{-1}$ this is equivalent to
\begin{equation} \label{eqn: twisted diff result2}
2(11)^{-1}\delta'_1(11)+
2(11)^{-1}\overline\gamma'(11) =
-2\varepsilon'_1 -2\overline\gamma'+
(k-2)((21)^{-1}H(21)-e_1^{-1}He_1).
\end{equation}

Finally we return to the equation $\Gamma=0$.
Observe that every term of $2\Gamma$ appears on the left-hand side of 
one of the equations
\eqref{eqn: twisted sum}, \eqref{eqn: twisted sum2},
\eqref{eqn: twisted diff result} or \eqref{eqn: twisted diff result2}
(where we set $i=1$ in \eqref{eqn: twisted sum},\eqref{eqn: twisted
sum2}) except for the terms $2\varepsilon'_1-2\delta'_1$.
Upon substituting the right-hand sides of 
\eqref{eqn: twisted sum}-\eqref{eqn: twisted diff result2} into
$2\Gamma$ most terms cancel and we are left with:
$$2\Gamma = 4\varepsilon'_1-4\delta'_1$$
and so $\varepsilon'_1=\delta'_1$.
Firstly taking the $\lt$ component of this implies $H=0$ 
(and so $\delta'_i=0=\varepsilon'_i$ for $i>1$) and secondly 
$\varepsilon'_1$ and $\delta'_1$ are in opposite Borels with zero
$\lt$ component and so must both be zero.
Now returning to equation \eqref{eqn: twisted diff result} we see
$$(11)\overline\gamma'(11)^{-1} = -\overline\gamma'$$
which says precisely that $\Ad_{\mu(p)}\gamma' = -\gamma'$
since $\mu(p) = C^{-1}(11)C$ in the notation we are using.
\epfms

\end{section}

\renewcommand{\baselinestretch}{1}              %
\normalsize
\bibliographystyle{amsplain}    \label{biby}
\bibliography{../thesis/syr}    
\end{document}